\documentstyle[psfig,12pt]{amsart}
\newcommand{\skipthistext}[1]{}
\newcommand{\remark}[1]{}
\renewcommand{\fbox}{}
\newcommand{\subtit}[1]{\vspace{0mm}

\noindent {\bf\large {#1}}

}
\newcommand{\rel}{\text{$\operatorname{rel}$}}

\newcommand{\Nd}{\text{$\operatorname{Nd}$}}
\newcommand{\card}{\text{$\operatorname{card}$}}
\renewcommand{\int}[1]{\stackrel{\circ}{{#1}}}
\newcommand{\diffeo}{\cong}
\newcommand{\defeq}{\stackrel{\text{def}}{=}}
\newcommand{\bcs}{\natural}
\renewcommand{\#}{\sharp}
\newenvironment{lemma}{\vspace{5mm}

\noindent {\bf Lemma.}\it}{\rm\vspace{5mm}

}

\newenvironment{theorem}{\vspace{5mm}

\noindent {\bf Theorem.}\it}{\rm\vspace{5mm}

}

\newenvironment{fact}{\vspace{2mm}

\noindent {\bf Fact.}\it}{\rm\vspace{0mm}

}

\newenvironment{definition}{\vspace{2mm}

\noindent {\bf Definition.}\it}{\rm\vspace{5mm}

}
\newcommand{\SPMXY}{
\vspace{5mm}

\centerline{
\fbox{
\begin{picture}(200,130)
    \put(0,0)       {\fbox{\psfig{figure=SPMXY1.ps}}}
    \put(52,92)     {$Y$}
    \put(12,05)     {$X$}
    \put(25,95)     {$S_*$}
    \put(64,110)    {$P_*$}
    \put(50,38)     {$M$}
    \put(175,-10)   {\makebox(0,0)[t]{\bf\large Figure 1.}}
\end{picture}
\begin{picture}(200,130)
    \put(5,0)       {\fbox{\psfig{figure=SPMXY2.ps}}}
    \put(20,92)     {$S_*$}
    \put(150,86)    {$P_*$}
    \put(32,50)     {$Y$}
    \put(113,50)    {$X$}
    \put(75,95)     {$M$}
\end{picture}
}
}
\vspace{10mm}

\noindent
}

\newcommand{\Trick}{
\vspace{5mm}

\centerline{
\fbox{
\begin{picture}(250,100)
    \put(5,0)      {\fbox{\psfig{figure=Trick.ps}}}
    \put(35,70)    {$M$}
    \put(35,30)    {$W_1$}
    \put(80,70)    {$W_1$}
    \put(80,30)    {$W_2$}
    \put(145,70)   {$W_1$}
    \put(145,30)   {$W_2$}
    \put(200,70)   {$M$}
    \put(200,30)   {$W_1$}
    \put(55,95)    {$M_1$}
    \put(180,95)   {$M_2$}
    \put(116,52)   {$S^3$}
    \put(114,-2)   {$\Sigma\sharp\Sigma$}
    \put(120,-10)  {\makebox(0,0)[t]{\bf\large Figure 2.}}
\end{picture}
}
}
\vspace{10mm}

\noindent
}

\newcommand{\SxSFinger}{
\vspace{5mm}

\centerline{
\fbox{
\begin{picture}(200,100)
    \put(3,0)       {\fbox{\psfig{figure=S2xS2.ps}}}
    \put(2,-2)      {$0$}
    \put(55,-2)     {$0$}
    \put(75,-2)     {$0$}
    \put(130,-2)    {$0$}
    \put(85,-10)    {\makebox(0,0)[t]{\bf\large Figure 3.}}
\end{picture}
\begin{picture}(200,100)
    \put(3,0)       {\fbox{\psfig{figure=FingerMove.ps}}}
    \put(0,70)      {\makebox(0,0)[rb]{\small  One of}}
    \put(0,60)      {\makebox(0,0)[rb]{\small  $S$'s}}
    \put(20,70)     {\makebox(0,0)[lb]{\small  one of}}
    \put(20,60)     {\makebox(0,0)[lb]{\small  $P$'s or $M$'s}}
    \put(175,90)    {\makebox(0,0)[lb]{\small  Accessory}}
    \put(175,80)    {\makebox(0,0)[lb]{\small  circle}}
    \put(175,45)    {\makebox(0,0)[lb]{\small  Whitney}}
    \put(175,35)    {\makebox(0,0)[lb]{\small  circle}}
    \put(105,25)    {\makebox(0,0)[rb]{\tiny   Pair of intersection}}
    \put(105,15)    {\makebox(0,0)[rb]{\tiny   points of opposite sign}}
    \put(100,-10)   {\makebox(0,0)[t]{\bf\large Figure 4.}}
\end{picture}
}
}
\vspace{10mm}

\noindent
}

\newcommand{\CuspKirby}{
\vspace{5mm}

\centerline{
\fbox{
\begin{picture}(200,170)
    \put(0,0)      {\fbox{\psfig{figure=CuspMove.ps}}}
    \put(-3,60)    {\makebox(0,0)[rb]{\small One of}}
    \put(-3,50)    {\makebox(0,0)[rb]{\small $S$'s}}
    \put(70,28)    {\makebox(0,0)[rb]{\tiny Selfintersection}}
    \put(70,18)    {\makebox(0,0)[rb]{\tiny point}}
    \put(173,45)   {\makebox(0,0)[lb]{\tiny Whitney}}
    \put(173,35)   {\makebox(0,0)[lb]{\tiny circle}}
    \put(80,-10)  {\makebox(0,0)[t]{\bf\large Figure 5.}}
\end{picture}
}
\fbox{
\begin{picture}(200,180)
    \put(0,0)      {\fbox{\psfig{figure=Kirby.ps}}}
    \put(-10,160)  {$P$}
    \put(135,128)  {$S$}
    \put(135,95)   {$M$}
    \put(20,25)    {$X$}
    \put(65,25)    {$Y$}
    \put(135,37)   {\makebox(0,0)[lb]{\tiny Attaching circles of}}
    \put(135,27)   {\makebox(0,0)[lb]{\tiny extra 2-handles}}
    \put(135,17)   {\makebox(0,0)[lb]{\tiny up to homotopy in $V_2$}}
    \put(150,152)  {\makebox(0,0)[lb]{\tiny Attaching circles of}}
    \put(150,142)  {\makebox(0,0)[lb]{\tiny extra 2-handles}}
    \put(150,132)  {\makebox(0,0)[lb]{\tiny up to homotopy in $\partial V_2$}}
    \put(110,-10)  {\makebox(0,0)[t]{\bf\large Figure 6.}}
\end{picture}
}
}
\vspace{10mm}

\noindent
}

\newsavebox{\hsmall}
\savebox{\hsmall}{\bf\large h}
\newcommand{\h}{\usebox{\hsmall}}
\title{A Decomposition of Smooth Simply-connected \h-cobordant 4-manifolds.}
\author{R. Matveyev}
\address{Address: Department of Mathematics\\ Michigan State University\\
E. Lansing, MI 48824.\\ matveyev@@math.msu.edu}
\thanks{This paper is a part of author's Ph.D. thesis.}
\date{}

\begin{document}
\maketitle
\vspace{15mm}
\subtit{Introduction and  the Statement.}

In \cite{A} S. Akbulut obtained an example of the exotic manifold cutting off
the contractible submanifold from the standard manifold and regluing it
via nontrivial involution of the boundary.

In these notes we give a proof of a decomposition theorem stated below.
It generalizes the example of Akbulut.

Another proof of the theorem was
independently obtained by C.L. Curtis, M.H. Freedman, W.C. Hsiang, and
 R. Stong in \cite{CHFS}.

The author would like to thank S. Akbulut for many useful discussions,
constant strong support and for
bringing \cite{CH0,CH1} to his attention.

Throughout these notes all maps and manifolds are smooth and immersions are
in general position (or do their best if they have to obey some extra
conditions).
We also make the convention that if a star appears in place of a subindex,
we consider a union of all objects in the family, where the index
substituted by the star runs over its range.
For example, $D_*\defeq\bigcup_i D_i$.
\begin{theorem}
Let
$U$ be a smooth, 5-dimensional, simply-connected h-cobordism with $\partial
U=M_1\sqcup(-M_2)$.
Let $f:M_1\rightarrow M_2$ be the homotopy equivalence induced
by $U$.
\begin{enumerate}
\item There are decompositions
\[
M_1=M\#_\Sigma W_1 , \hspace{5mm}
M_2=M\#_\Sigma W_2
\]
such that $in_{2*}^{-1} \circ in_{1*} = f_* : H_2(M_1) \rightarrow H_2(M_2)$.
Here  $in_{2*}$, $in_{1*}$ are the maps induced in the second homology
by embeddings of $M$ into $M_1$ and $M_2$ respectively, and $W_1$, $W_2$ are
smooth, compact, contractible 4-manifolds, and
$\Sigma=\partial W_1=\partial W_2=\partial M$.
\item These decompositions may be chosen so that $W_1$ is diffeomorphic to
$W_2$.
\end{enumerate}
\end{theorem}
In fact, it can be seen from the proof that the whole cobordism can be
decomposed
into two subcobordisms, each is a product cobordism and one is
diffeomorphic to $D^4$ (as a smooth manifold, without any additional
structure).

We will also need the following
\begin{definition}
We say that two collections $\{S_i\}_{i=1}^n$, $\{P_i\}_{i=1}^n$ of oriented
2-spheres  immersed in oriented 4-manifold are algebraically dual if
$\langle[S_i],[P_j]\rangle=\delta_{ij}$.
Here, $[S]$ is a homology class of immersed sphere $S$ and
$\langle \cdot\, , \cdot \rangle$ is the intersection form in the second
homology
of 4-manifold.

They are geometrically dual if they are algebraically dual and, moreover,
$\card(S_i\cap P_j)=\delta_{ij}$.
\end{definition}

\subtit{Proof of Theorem.}
First, observe that $U$ has  a handlebody with no 1- and 4-handles.
Let $N$ be the middle level of $U$ between 2- and 3-handles.
Then we have two diffeomorphisms

\centerline{$\phi_1:M_1\#S_{11}^{2}\times S_{12}^{2}\#\ldots\#S_{n1}^{2}
\times S_{n2}^{2}\rightarrow N$}
\centerline{$\phi_2:M_2\#S_{11}^{2}\times S_{12}^{2}\#\ldots\# S_{n1}^{2}
\times S_{n2}^{2}\rightarrow N$.}

By enriching the handlebody of $U$ by 2-3 canceling pairs of handles and
choosing
$\phi_1$ and $\phi_2$, we can assume that
$(\phi_2^{-1}\circ\phi_1)_*|_{H_2(M_1)}=f_*$ and
$(\phi_2^{-1}\circ\phi_1)_* [S_{ij}^2]=[S_{ij}^2]$,
$i=1,\ldots,n$, $j=1,2$.
Then the two embeddings
$c_1:\bigsqcup S_{i1}^2\vee S_{i2}^2\rightarrow N$,
$c_2:\bigsqcup S_{i1}^2\vee S_{i2}^2\rightarrow N$,
representing cores of products of spheres in decompositions $\phi_1$ and
$\phi_2$, are homotopic.
We have two algebraically dual collections of embedded 2-spheres
$\{S_i=c_1(S_{i1}^2)\subset N\}_{i=1}^n$, $\{P_i=c_2(S_{i2}^2)\subset
N\}_{i=1}^n$,  and
surgery along   $\{P_i\}$ gives $M_1$; along  $\{S_i\}$ gives $M_2$.
Put $V_0=\overline{\Nd_N(S_*\cup P_*)}$, the closed regular neighborhood of
$S_*\cup P_*$ in $N$.

Here we give a sketch of the rest of the construction and than work out the
details.
\begin{itemize}
\item[{\bf 1.}] Manifold $V_0$ has a free fundamental group and its second
     homology are generated by classes of spheres $S_i$ and $P_i$,
     $i=1,\ldots ,n$.
     We enlarge $V_0$ to obtain the bigger manifold $V_1$ so that collections
     of spheres $\{S_i\}$ and $\{P_i\}$ have geometrically dual collections of
     immersed spheres inside of $V_1$  and the properties mentioned above are
     satisfied.
\item[{\bf 2.}] We add 1-handles and essential 2-handles to $V_1$ so that the
     fundamental group vanishes and   no  new two dimensional
     homology classes appear.
\item[{\bf 3.}] Surgery of the manifold obtained in step 2 along collections
     $\{S_i\}$ and  $\{P_i\}$ gives  two contractible submanifolds $W_1$ and
     $W_2$ of $M_1$ and $M_2$, respectively, and
     $\overline{M_1\setminus W_1}\diffeo \overline{M_2\setminus W_2}$.
\end{itemize}

The intersection points of embedded spheres $S_i$ and $P_j$ can be grouped in
pairs of points of opposite sign with one extra point of positive sign when
$i=j$, which we refer to as  wedge points.
Consider a disjoint collection of Whitney circles $l_1,l_2,\ldots,l_m$ in
$S_*\cup P_*$, one for each pair of points considered above.
Push these circles to the boundary of the neighborhood of $S_*\cup P_*$ in
$N$ and call this new circles $l'_1,l'_2,\ldots,l'_n$ as in.
Each of $l'_i$ is contractible in $N\setminus S_*$, as well as in
$N\setminus P_*$, since the latter two are simply-connected manifolds.
Thus, for each  $l_i$ we can find two immersed disks $D_i$ and $E_i$ with
$\partial D_i=\partial E_i=l_i$, coinciding along the collar of the
boundary, and satisfying the following condition:
\begin{equation}
\int{D_*}\cap  S_*=\O,\hspace{5mm}
\int{E_*}\cap  P_*=\O .\label{cond:disj}
\end{equation}

\skipthistext{
Applying Whitney tricks to $P_*$ along $\{D_i\}$ gives as a collection of
immersed spheres geometrically dual to $\{S_i\}$. Same way we can obtain
collection of spheres geometrically dual to $\{P_i\}$ applying Whitney
trick to $S_*$ along $E$'s.
But we can not add both $D_*$ and $E_*$ to our manifold because we don't
want to create new two dimensional homology classes.
}

Now we shall show that it is possible to find disks $D_i'$ and $E_i'$
so that they are homotopic \rel(collar of the boundary) and the condition
similar to (\ref{cond:disj}) is satisfied.

The union of $D_i$ and $E_i$ with appropriate orientations gives us a class
$[D_i\cup E_i]$ in the homology of $N$, which splits into the direct sum
$H_2(N)=H_2(M_1)\oplus<[S_i],[P_i]>$.
So, we can write
$[D_i\cup E_i]=a+\sum{}{}\beta_k  [S_k]+\sum{}{}\gamma_k [P_k]$, where
$a\in H_2(M_1)$.
Using the fact that $\pi_1(M_1)=\{1\}$ and the decomposition $\phi_1$ we can
realize
class $a$, considered as a class in $H_2(N)$, by an immersed sphere $A$ in $N$
disjoint from  $S_*$.
Classes $\sum{}{}\beta_k  [S_k]$ and $\sum{}{}\gamma_k [P_k]$ can be
realized by the immersed spheres $B$ and $C$ in $N$ disjoint from $S_*$ and
$P_*$, respectively.
Now taking  the connected sums $D'_i=D_i\#(-A)\#(-B)$, $E'_i=E_i\#(-C)$
(`$-$' here denotes reversal of orientation) ambiently along carefully
chosen paths, we obtain discs $D'_i$ and $E'_i$ satisfying property similar
to (\ref{cond:disj}) and homotopic \rel(collar of the boundary).

Consider homotopy  $F_t:\bigsqcup D_i^2\rightarrow N$
\rel(collar of the boundary), where $F_0(\bigsqcup D_i^2)=D_*'$ and
$F_1(\bigsqcup D_i^2)=E_*'$.
It can be also viewed as a homotopy of the union of the disks and spheres
$S_*$,
$P_*$, where the spheres stay fixed during the homotopy.
\begin{lemma}
Homotopy $F$ can be perturbed, fixing ends, $S_*\cup P_*$ and the collar of
the
boundary of disks in the disks, to homotopy $F'$ satisfying the following
property: $F|_{[0,\frac{1}{2}]}$ can be decomposed in a sequence
of simple homotopies each of which is either a cusp or finger move;
analogously,  $F|_{[\frac{1}{2},1]}$ can be decomposed in a sequence of
inverse cusp moves and Whitney tricks.
\end{lemma}

For definitions of  cusp, finger move, inverse cusp move and Whitney trick
see~\cite{FQ,K,GM}.

Proof of  Lemma is given in the next section.

Assume, now, that $F$ has the property provided by Lemma above.
Consider $M_*=\bigcup_i M_i=F_{\frac{1}{2}}(\bigsqcup D^2_i)$.
Each of $M_i$ may intersect both $S_*$ and $P_*$, but it is homotopic
\rel(collar of the boundary) to the disk with interior disjoint from
$S_*$, via a homotopy increasing number of intersection points.
This means that the intersection points of $\int{M_*}$ and $S_*$
can be grouped
in pairs so that for each pair there is an embedded Whitney disk $X_k$
with interior disjoint from the rest of the picture.
The same argument shows that there are disks $Y_k$ for intersection points of
$M_*$ and $P_*$, see figure 1.
\SPMXY

Using embedded disks from collection $\{X_k\}$ we can push all disks $M_i$
off $S_*$, producing the new disks $M_i'$. Each $M_i'$ has an interior
disjoint
from $S_*$ and its boundary is a Whitney circle for a pair of  points of
opposite
sign in $S_*\cap P_*$.
Applying an immersed Whitney trick to $P_*$ we eliminate all the
intersections with  $S_*$, except those  required by algebraic conditions
and produce the collection $\{S_i^\perp\}$ geometrically dual to
$\{S_i\}$.
Note that, since disks $M_i'$ are, in general, immersed and may have
non-trivial relative normal bundle in $N$, we may have to introduce
(self-)intersection
points to spheres $\{S_i^\perp\}$ during this process.
In the same way, using disks from collection $\{Y_i\}$ we can obtain the
immersed Whitney disks $M_i''$ disjoint from $P_*$.
They allow us to create a collection $\{P_i^\perp\}$ of immersed
spheres dual to $\{P_i\}$.
Note, that all described homotopies can be performed in a regular
neighborhood of
$M_*\cup X_*\cup Y_*$.

Thus, one can see that the manifold $V_1=V_0 \cup \overline{\Nd_N(M_*\cup
X_*\cup Y_*)}$
has the following properties:
\begin{enumerate}
\item $H_2(V_1)$ is generated by homology classes of embedded spheres
      $S_i$, $P_i$.
\item The collection of spheres $\{S_i\}$ has the geometrically dual
collection
      $\{S_i^\perp\}$ of immersed spheres in $V_1$.
      Similarly, collection of spheres $\{P_i\}$ has geometrically dual
      collection $\{P_i^\perp\}$.
\item $\pi_1(V_1)$ is a free group.
\end{enumerate}

Consider a handlebody of $N$ starting from $V_1$.
Put $V_2=V_1\cup(\text{union of 1-handles})$.
$V_2$ still satisfies properties similar to (1), (2), (3) for $V_1$.

Let $g_1,\ldots ,g_l$ be free generators of $\pi_1(V_2)$.
If we fix paths from a basepoint to attaching spheres of 2-handles then
they represent elements of $\pi_1(\partial V_2)$, say $h_1,\ldots ,h_L$.
Since $V_2\cup(\text{all 2-handles})$ is a simply connected manifold, each
$g_i$ has a lift $\tilde{g}_i$ in $\pi_1(\partial V_2)$, such that it
belongs to the normal subgroup of $\pi_1(\partial V_2)$ generated by
$h_1,\ldots ,h_L$.
In other words
$\tilde{g}_i=(\alpha_1 h_{i_1}^{\pm 1}\alpha_1^{-1})(\alpha_2 h_{i_2}^{\pm
1}\alpha_2^{-1})\ldots (\alpha_k h_{i_k}^{\pm 1}\alpha_k^{-1})$.
Duplicating  2-handles, equipping them with appropriate
orientations and choosing paths joining attaching circles of 2-handles
to the basepoint
we may write $\tilde{g}_i=h_1' h_2'\ldots h_k'$, where
$h_1', h_2',\ldots, h_k'$ are elements of $\pi_1(\partial V_2)$ obtained
from attaching spheres of new handles with new paths to the basepoint.
Now, handles  in the decomposition of $\tilde{g_i}$ are all distinct and
we can slide the first handle over all others along paths joining
feet of the handles to the basepoint.
The attaching circle of the resulting
handle, call it $G_i$, is freely homotopic to $g_i$ in $V_2$.
Adding  such $G_i$'s to $V_2$ for each generator of $\pi_1(V_2)$ we obtain
the simply-connected manifold $V_3$.
Since $g_1,\ldots ,g_l$ are free
generators of the fundamental group, we do not create any additional
two-dimensional homology classes.

Surgery of $V_3$ along collections of embedded spheres $\{S_i\}$ and
$\{P_i\}$
gives two contractible sub-manifolds $W_1$ and $W_2$ of $M_1$ and $M_2$,
respectively.
Put
$M\defeq \overline{N\setminus V_3}\diffeo \overline{M_1\setminus W_1}
\diffeo \overline{M_2\setminus W_2}$.
So we have the decompositions:
\[
M_1=M\#_\Sigma W_1 , \hspace{5mm}
M_2=M\#_\Sigma W_2
\]
The property of induced maps in the homology stated in  Theorem is obvious.
The first part of  Theorem is proved.

\

Proof of the second part is based on the
\begin{fact}
If $W_1$, $W_2$ are homotopy balls built in the proof of the first part
of  Theorem and $S^4$ is a 4-dimensional sphere with standard smooth
structure, then
\[ W_1\#_\Sigma W_1 \diffeo S^4, \hspace{8mm} W_1\#_\Sigma W_2 \diffeo S^4. \]
\end{fact}
Denote the boundary connected sum by `$\bcs$'.
Then we can write
$$M_1=M\#_\Sigma W_1 \diffeo (M\#_\Sigma W_1)\#(W_1\#_\Sigma W_2)\diffeo
 (M\bcs W_1)\#_{\Sigma\#\Sigma}(W_1\bcs W_2)$$
and
$$M_2=M\#_\Sigma W_2 \diffeo (M\#_\Sigma W_2)\#(W_1\#_\Sigma W_1)\diffeo
 (M\bcs W_1)\#_{\Sigma\#\Sigma}(W_2\bcs W_1).$$
See figure 2.
\Trick

In order to proof Fact we will use the art of Kirby calculus.

First, we  build a handlebody of $V_3$.

Consider $V'_0$, a closed neighborhood of $S_*\cup P_*\cup(\text{arcs
joining wedge points in
$S_i\cap P_i$ and  $S_{i+1}\cap P_{i+1}$)}$.
If there were no intersections between $S_*$ and $P_*$, except a wedge point
for each pair $S_i$, $P_i$, then the handlebody would look like as shown on
figure 3.
\SxSFinger

Introducing a pair of intersection points of opposite signs corresponds
to the move in Kirby calculus
shown in figure 4.

Here, we  introduce two 1-handles, one corresponding to a Whitney circle,
another to an accessory circle of a newly introduced pair of intersections.
The handlebody of $V'_0$ is obtained by applying several moves shown in
figure 4 to the picture in figure 3.
Again, if disks $M_i$ were embedded and disjoint from $S_*\cup P_*$, then
to obtain the handlebody of $V_0\cup \overline{\Nd(M_*)}$ one has to attach
2-handles to Whitney circles for each pair of intersection points of $S_*$
and
$P_*$.
Then we have to introduce intersections between $M_*$ and $S_*\cup P_*$ and
self-intersections of $M_*$.
Corresponding moves are shown on figures 4 and 5.
\CuspKirby

Introducing a self-intersection corresponds to adding one 1-handle to the
picture.
Addition of disks from collections $\{X_i\}$, $\{Y_i\}$ corresponds to
attaching 2-handles to the circles linked once to the 1-handles
corresponding
to the Whitney circles,
and unlinked from other 1-handles.
They may link each other according to
the fact that boundaries of $X$'s and $Y$'s are not necessarily disjoint.
This phenomenon is illustrated on figure 2.
Then, as in previous steps, add 1-handles to the picture reflecting
intersections of $X_*$ and $Y_*$, as in figure 6.

To obtain the handlebody of $V_2$ we have to attach several  1-handles.
The link on figure 6 shows all possible phenomena which can occur.
Thus, the handlebody of $V_2$ has 1-handles of six types coming from
\begin{itemize}
\item[1.] Whitney circles of intersections of $S_*$ and $P_*$;
\item[2.] accessory circles of intersections of $S_*$ and $P_*$;
\item[3.] Whitney circles of intersections of $M_*$ and $S_*\cup P_*$;
\item[4.] accessory circles of intersections of $M_*$ and $S_*\cup P_*$;
\item[5.] Whitney and accessory circles of self-intersections of $M_*$;
\item[6.] Extra 1-handles.
\end{itemize}
For each 1-handle of type 1 and 3 there is a 2-handle attached to the circle
linked to this 1-handle geometrically once and algebraically zero times to
other
1-handles.

The attaching circles $h_1,\ldots, h_l$ of other 2-handles are homotopic to
free generators of $\pi_1(V_2)$.
As generators we may choose the cores
$g_1,\ldots, g_l$ of 1- handles of type 2, 4, 5 and 6.

Consider the homotopy $F:S^1\times I\rightarrow V_2$, $F(\cdot,0)=g_i$,
$F(\cdot,1)=h_i$.
We can make it disjoint from $S_*\cup M^*\cup X_*\cup Y_*$.
First, intersections of the image of the homotopy with $X_*$ and $Y_*$ can be
turned into intersections with $M_*$ by pushing them toward the boundary of
$X_*\cup Y_*$. Then, in the same way, we can avoid intersections with $M_*$
by  cost of new intersections with $P_*$.
Finally, intersections with
$S_*$ can be removed using the  geometrically dual collection of immersed
spheres $\{S_i^\perp\}$.
Call this new homotopy $F'$ and let $\{x_1,x_2,\dots,x_k\}={F'}^{-1}(P_*)$.
If we remove a small  neighborhood of a union of disjoint arcs joining each
$x_i$
to the point on $S^1\times \{0\}$ from $S^1\times I$, then restriction of
$F'$ on this set is a homotopy of $h_i$ to the curve ${g'}_i$ which is a band
connected sum  of $g_i$ and meridians of $P_*$.

This homotopy is disjoint from
$S_*\cup P_*\cup M_*\cup X_*\cup Y_*$ and can be pushed to the boundary of
$V_2$.

We obtain the handlebody of $W_1$ by attaching 2-handles to $h_1,\ldots, h_l$
and surgering $S_*$, which corresponds to putting  dots on the circles
representing $S_*$ in Kirby calculus.

Here we summarize all the information about the handlebody of $W_2$.
For 2-handles we use same notation as for the corresponding object in
construction
of $W_1$ above.
\newlength{\gap}
\setlength{\gap}{1mm}
\newlength{\leftentry}
\setlength{\leftentry}{75mm}
\newlength{\rightentry}
\setlength{\rightentry}{75mm}
%
%
\newcommand{\tit}{\parbox{\leftentry}{\centerline{1-handles}}}
%
%
\newcommand{\tittit}{\parbox{\rightentry}{\centerline{2-handles}}}
%
%
\newcommand{\SS}{\parbox{\leftentry}{Surgery of $S_*$.}}
%
%
\renewcommand{\P}{\parbox{\rightentry}{\vspace{\gap}$P_*$; attaching
circles
go geometrically
                                         once through corresponding 1-handles
                                         $S_*$.\vspace{\gap}}}

%
%
\newcommand{\SnP}{\parbox{\leftentry}{Whitney circles of intersections of
$S_*$
                                      and $P_*$.}}
%
%
\newcommand{\M}{\parbox{\rightentry}{\vspace{\gap}$M_*$; attaching circles go
                                 geometrically
                                 once through corresponding 1-handles from
                                 the left entry of the table.\vspace{\gap}}}
%
%
\newcommand{\MnSP}{\parbox{\leftentry}{Whitney circles of intersections of
                                    $M_*$ and $S_*\cup P_*$.}}
%
%
\newcommand{\XY}{\parbox{\rightentry}{\vspace{\gap}$X_*$, $Y_*$;
attaching circles go geometrically
                                  once through corresponding 1-handles from
                                  the left entry of the table.\vspace{\gap}}}
%
%
\newcommand{\ExtraOne}{\parbox{\leftentry}{Accessory circles of all
                                        intersections, Whitney circles of
                                        (self-)intersections of $M_*$, Whitney
                                        circles of intersections of $X_*$
                                        and $Y_*$, extra 1-handles.}}
%
%
\newcommand{\ExtraTwo}{\parbox{\rightentry}{\vspace{\gap}Extra 2-handles $H_*$;
                                        attaching circles
                                        go geometrically once through
                                        corresponding 1-handles from  the
                                        left entry of the table and homotopic
                                        to band connected sum of 1-handles
                                        with meridians of $P_*$.\vspace{\gap}}}
\vspace{12pt}
\noindent
\begin{tabular}{lll}
  & \tit & \tittit \\
&&\\
1.& \SS                     & \P                      \\
&&\\
2.& \SnP                    & \M                      \\
&&\\
3.& \MnSP                   & \XY                     \\
&&\\
4.& \ExtraOne               & \ExtraTwo               \\
\end{tabular}
\vspace{12pt}

Taking a double of $W_1$ corresponds to attaching zero-framed 2-handles to
the meridians of existent 2-handles, and as many 3-handles as there
are 1-handles in the handlebody of $W_1$.

Remember that attaching circles $h_i$ of 2-handles $H_i$ are homotopic to
$g_i'$.
Sliding $H_*$ over dual handles $H_*^*$ we can obtain handles $H_*'$ attached
to
$g_*'$.
Now ${g'}_1,\ldots, {g'}_l$ may be linked
to the meridians of handles $P_*$, $M_*$, $X_*$ and $Y_*$,
but situation  can be improved by sliding handles dual to $P_*$, $M_*$, $X_*$
and $Y_*$ over  handles dual to $H'_*$.
Further we may slide $H'_*$ over handles dual to $P_*$ to obtain handles
attached to meridians of 1-handles from the fourth row of the table, so they
can be
cancelled.
The framings of $H'_*$ result in a twist of attaching circles of remaining
2-handles.

We undo this twist using dual handles.

With the next step we unlink handles $X_*$ and $Y_*$ from each other
using dual 2-handles and
cancel them with 1-handles of type 3. Now $M_i$ are attached to meridians
of 1-handles of type 1 and also can be cancelled. After cancelling $P_*$
with 1-handles from surgery of $S_*$ we end up with unknotted,
unlinked, zero framed handles and the same number of 3-handles, which is
clearly $S^4$.

We need only small changes in our construction to prove the second
diffeomorphism in our Fact.
The handlebody of $W_2$ differs from $W_1$ by putting  dots on the
attaching circles of $P_*$ rather than $S_*$.
Thus, the handlebody of $W_1\#_\Sigma W_2$ is obtained from $W_1$
by attaching 2-handles to meridians of 1-handles, coming from surgered
$S_i$'s, and to meridians of all 2-handles, except $P_i$'s.
Using the same trick we can make the homotopy of $h_*$ to $g_*$ disjoint
from $P_*\cup M_*\cup X_* \cup Y_*$, rather than $S_*\cup M_*\cup X_* \cup
Y_*$.
Applying the same procedure to $H_*$ gives us 2-handles $G''_*$ attached
to the
band connected sum of  meridians of 1-handles corresponding to generators of
$\pi_1(V_2)$ and meridians of $S_*$.
But now 1-handles corresponding
$S_*$ have dual 2-handles, so sliding
$H_*$ over them gives handles dual to 1-handles generating $\pi_1(V_2)$.
So we may apply the same procedure to simplify the
handlebody and end up with $S^4$.
This finishes the proof of Fact and the second part of Theorem.

\subtit{Proof of Lemma.}
It is a simple consequence of singularity theory that a homotopy
of a surface in 4-manifold can be decomposed (after small
perturbation) in a sequence of finger moves, cusps, Whitney tricks
and inverse cusp moves.
In our case we have to consider, in addition, a finger move along a path
in one of the disks joining (self-)intersection of the disks with a point on
its
boundary.
The inverse of this move is a Whitney trick with a Whitney circle intersecting
the boundary of the disk.
We have to show that it is possible to reorder these simple homotopies so
that those increasing number of intersections come first.
This is obvious after the following consideration: cusp birth happens
in a small neighborhood of the point in the surface and we can
assume that part of the surface in this neighborhood stays fixed
during the part of homotopy preceding this cusp birth, so we can push
this cusp
up to the beginning of the homotopy.
A finger move can be localized in the neighborhood of an embedded arc
joining two points in the surface, thus we may apply the same argument
and this finishes the  proof of Lemma.

\end{document}